\documentclass[dvips]{arxbj}
\usepackage{upgreek}
\usepackage{graphicx}

\aid{0}
\volume{16}
\issue{3}
\pubyear{2010}
\firstpage{780}
\lastpage{797}
\doi{10.3150/09-BEJ226}

\makeatletter
\renewcommand{\citep}[1]{(\citet{#1})}

\newcommand{\1}{\mathbf{1}}
\newcommand{\DD}{\mathrm{d}}
\newcommand{\starexp}{\exp^*}
\newcommand{\B}{\mathcal{B}}
\newcommand{\CC}{\mathbb{C}}
\newcommand{\Cov}{\operatorname{Cov}}
\newcommand{\D}{\,\mathrm{d}}
\newcommand{\diag}{\operatorname{diag}}
\newcommand{\EE}{\mathbb{E}}
\newcommand{\RR}{\mathbb{R}}
\newcommand{\N}{\mathcal{N}}
\newcommand{\NN}{\mathbb{N}}
\newcommand{\Var}{\operatorname{Var}}

\newremark{remark}[thm]{Remark}
\newtheorem{proposition}[thm]{Proposition}
\newtheorem{lemma}[thm]{Lemma}
\newremark{example}[thm]{Example}

\makeatother

\begin{document}
\begin{frontmatter}

\title{Some covariance models based on normal scale mixtures}
\runtitle{Covariance models based on normal scale mixtures}

\begin{aug}
\author{\fnms{Martin} \snm{Schlather}\ead[label=e1]{schlather@math.uni-goettingen.de}}
\runauthor{M. Schlather}
\address{Institut f\"{u}r Mathematische Stochastik \& Zentrum f\"{u}r
Statistik,
Universit\"{a}t G\"{o}ttingen, Goldschmidtstr. 7, D-37077 G\"
{o}ttingen, Germany. \printead{e1}}
\end{aug}

\received{\smonth{11} \syear{2008}}
\revised{\smonth{5} \syear{2009}}

%
\begin{abstract}
Modelling spatio-temporal processes
has become an important issue in current research.
Since Gaussian processes are essentially determined
by their second order structure, broad classes of
covariance functions are of interest.
Here, a new class is described that
merges and generalizes various models presented in the literature,
in particular models in
Gneiting (\textit{J. Amer. Statist. Assoc.} \textbf{97} (2002)
590--600) and Stein (Nonstationary spatial covariance functions (2005)
Univ. Chicago).
Furthermore, new models and a multivariate extension are introduced.
\end{abstract}

%

\end{frontmatter}

\section{Introduction}
Spatio-temporal modelling is an important task in many disciplines of
the natural
sciences, geosciences, and engineering. Hence,
the development of models for spatio-temporal
correlation structure
is of particular interest. The lively activity in this field of research
has become apparent through
various recent reviews of known classes of
spatio-temporal covariance functions (\citet{GGG07}, \citet{MPG08},
\citet{ma08}).
To categorise these classes, different
aspects have been considered. \citet{GGG07} distinguish between
the properties of covariance functions, such as motion invariance, separability,
full symmetry, or conformity with
Taylor's hypothesis. Another classification is based on the
construction principles (\citet{ma08}), such as
spectral methods \citep{stein05},
multiplicative mixture models \citep{ma02},
additive models \citep{ma05c},
turning bands upgrade \citep{KCHS04},
derivatives and integrals \citep{ma05b}, and
Gneiting's (\citeyear{gneitingnonseparable}) approach,
see also \citet{stein05c} and \citet{ma03}.

Surprisingly, some rather different approaches to the construction of
spatial and spatio-temporal covariance models
can be subsumed in a unique class of normal scale mixtures,
which is a generalization of Gneiting's
(\citeyear{gneitingnonseparable}) class.
As its construction is based on cross covariance functions,
Section \ref{sec:background} illustrates some of the properties of cross
covariance
functions and cross variograms.
In Section~\ref{sec:gneiting}, Gneiting's class itself is
generalized. Section~\ref{sec:random}
introduces two new classes of spatio-temporal models.
Section \ref{sec:multivariate} presents an extension
to multivariate models.
In addition to the two-dimensional realisations illustrated below,
three-dimensional realisations
are available in the form of films at the following website:
\href{http://www.stochastik.math.uni-goettingen.de/data/bernoulli10/}{www.stochastik.math.uni-goettingen.de/data/ bernoulli10/}.


\section{Background: Cross covariance functions}
\label{sec:background}
Here we introduce some basic notions and
properties of cross covariance functions and cross variograms.
See \citet{wackernagel} for a geostatistical overview and
\citet{reisertburkhardt07} for some of the construction principles
of multivariate cross covariance functions in a general framework.

Let $Z(x)=(Z_1(x), \ldots, Z_m(x))$, $x \in\RR^d$, be a zero mean,
second order
$m$-variate, complex
valued random field in $\RR^d$, that is,
$\Var Z_j(x)$ exists and
$\EE Z_j(x) =0$
for all
$x\in\RR^d$ and $j=1,\ldots, m$.
Then, the cross covariance function $C\dvtx\RR^{2d}\rightarrow\CC
^{m\times m}$
is defined by
\[
C_{jk}(x, y) = \Cov(Z_j(x),Z_k(y)),\qquad x,y\in\RR^d,
j,k=1,\ldots,m
.
\]
Clearly $C(x,y) = \overline{C^\top(y,x)}$,
but $C(x,y)= \overline{C^\top(x,y)}$ is not valid
in general.
A function $C\dvtx\RR^{2d} \rightarrow\CC^{m\times m}$ with $C(x,y) =
\overline{C^\top(y,x)}$, $x,y\in\RR^d$, is
called positive definite
if for all $n\in\NN$, $x_1,\ldots, x_n\in\RR^d$ and
$a_1,\ldots,a_n\in\CC^m$,
\begin{equation}
\label{eq:posdef}
\sum_{p=1}^n\sum_{q=1}^n a_p^\top C(x_p, x_q) \bar a_q \ge0.
\end{equation}
It is called strictly
positive definite if strict inequality holds in (\ref{eq:posdef}) for
$(a_1,\ldots,a_n)\not= 0$ and pairwise distinct $x_1,\ldots, x_n$.
Accordingly, we name a Hermitian matrix $M\in\CC^{m\times m}$ positive
definite, if
$v^\top M \bar v\ge0$ for all $v\in\CC^m$, and strictly positive definite
if strict inequality holds for $v\not=0$.

As in the univariate case, we derive from Kolmogorov's existence
theorem that
a function
$C\dvtx\RR^{2d} \rightarrow\CC^{m\times m}$ with $C(x,y) = \overline
{C^\top(y,x)}$
is a positive
definite function
if and only if a (Gaussian) random field
exists with $C$ as cross covariance function.
Further, a
function $C \dvtx\RR^{2d} \rightarrow\RR^{m\times m}$ is a positive
definite function if and only if Equation~(\ref{eq:posdef}) holds for
any $a_1,\ldots,a_n\in\RR^m$.

The cross variogram $\gamma\dvtx\RR^{2d}\rightarrow\CC^{m\times
m}$, $\gamma=(\gamma_{jk})_{j,k=1,\ldots,m}$
is defined by
\[
\gamma_{jk}(x, y) = \tfrac12 \EE\bigl(
Z_j(x)-Z_j(y)\bigr)\overline{\bigl(Z_k(x)-Z_k(y)\bigr)},\qquad
x,y\in\RR^d, j,k=1,\ldots,m.
\]
If $Z$ has second order stationary
increments, then $\gamma(x,y)$ depends only on the distance
vector $h=x-y$, that is, $\gamma(x,y)=\tilde\gamma(h)$ for some
function $\tilde
\gamma\dvtx\RR^d \rightarrow\CC^{m\times m}$.
If in addition $Z$ is univariate,
then $\tilde\gamma$ is called a (semi-)variogram.
\citeauthor{schoenberg38b}'s (\citeyear{schoenberg38b}) theorem
states that a
function $\tilde\gamma\dvtx\RR^{d}\rightarrow\RR$ with $\tilde\gamma
(0)=0$ is
a variogram if and only if $\exp(-r\tilde\gamma)$ is a covariance function
for all $r>0$, see also \citet{GSS01}.
Let us now discuss
multivariate and non-stationary versions of this
statement.
To this end, we denote the componentwise multiplication of matrices by ``$*$'',
in particular,
\[
A^{*n} = (A_{jk}^n)_{jk}\qquad
\hbox{for } A= (A_{jk})_{jk}
.
\]
Further, $f^*(A)$ denotes the componentwise function evaluation, for example,
\[
\starexp(A) = (\exp(A_{jk}))_{jk}.
\]

\begin{thm}\label{thm:Mplus}
Let $C\dvtx\RR^{2d}\rightarrow\CC^{m\times m}$ and $E_{m\times m}$ be the
$m\times
m$ matrix whose components are all~$1$.
\begin{enumerate}
\item
The following three assertions are equivalent:
\textup{(i)} $C$ is a cross covariance function; \textup{(ii)}~$\exp^*(rC)-E_{m\times m}$ is a cross covariance function for all $r>0$;
\textup{(iii)}
$\sinh^*(rC)$ is a cross covariance function for all $r>0$.
\item
If $\exp^*(rC)$ is a cross covariance function for all $r>0$
then
\begin{equation}\label{eq:x0}
C^{(z)}(x,y) = C(z, z) -C(x, z) -C(z, y) + C(x,y)
\end{equation}
is a cross covariance function
for all $z\in\RR^d$.
If $m=1$ and (\ref{eq:x0}) holds for one $z\in\RR^d$, then
$\exp(rC)$ is a covariance function for all $r>0$.
\end{enumerate}
\end{thm}

\begin{pf}
Note that the componentwise product $C_1 * C_2$
of two $m$-variate
cross covariance functions $C_1$ and $C_2$ is again a cross covariance
function.
To see this, consider
the componentwise product of two independent random fields with
cross covariance functions $C_1$ and $C_2$.
In particular, $C(x,y)^{*n}$ and $r C(x,y)$, $r\ge0$, are cross covariance
functions.
Furthermore, the sum and the pointwise
limit of $m$-variate cross covariance functions
are cross covariance functions.
\begin{longlist}[2.]
\item[1.]
Both functions, $\exp(x)-1$
and $\sinh(x)$, have Taylor expansion on $\RR$
with positive coefficients only.
Hence, $\starexp(rC) -E_{m\times m}$ and $\sinh^*(rC)$
are cross covariance functions
if $C$ is a cross covariance function.
On the other hand, since the Taylor expansions
equal $x+\mathrm{o}(x)$ as $x\rightarrow0$, we have that
$(\starexp(rC)
-E_{m\times m})/r$ and $\sinh^*(rC)/r$ converge to $C$ as
$r\rightarrow0$
and $C$ must be a cross covariance function.
\item[2.]
The proof follows the lines in \citet{matheron72}.
Let $a_1,\ldots, a_n\in\CC^m$, $x_1,\ldots, x_n\in\RR^d$,
$a_0= -\sum_{p=1}^n a_p$ and $x_0=z$ for some $z\in\RR^d$. Then
\begin{eqnarray*}
0& \le&
\lim_{r\rightarrow0} \sum_{p=0}^n \sum_{q=0}^n
a_p^\top\frac{\starexp(rC(x_p,x_q))-
E_{m\times m}}{r}
\bar a_q
=
\sum_{p=0}^n \sum_{q=0}^n a_p^\top C(x_p, x_q)\bar a_q
\\
&=&
\sum_{p=1}^n \sum_{q=1}^n a_p^\top
[C(x_p, x_q) + C(z, z) - C(x_p, z) - C(z, x_q)
]\bar a_q.
\end{eqnarray*}
Conversely, assume that $m=1$ and Equation~(\ref{eq:x0}) holds.
Since $C_0(x,y) = f(x)\overline{f(y)}$
is a covariance function for any function
$f\dvtx\RR^d\rightarrow\CC$ (\citet{BTA}, Lemma 1)
part 1 of the theorem results in
\[
\exp(rC(x, y)) =
f(x) \overline{f(y)} \exp\bigl(rC(x, y) + rC(z, z) - rC(x, z) - rC(z, y)\bigr)
\]
being a positive definite function
for any $r>0$ and $f(x) = \exp(rC(x,z) - rC(z,z) /2)$.\qed
\end{longlist}\noqed
\end{pf}

\begin{remark}
If $m=1$, $C(x,y)= -\tilde\gamma(x-y)$ and $z=0$,
then $C^{(0)}$ in Equation~(\ref{eq:x0}) equals the covariance function of an
intrinsically
stationary random field $Z$ with $Z(0)=0$ almost surely, that is,
part 2 of Theorem \ref{thm:Mplus} yields
\citeauthor{schoenberg38b}'s (\citeyear{schoenberg38b}) theorem.
If $m>1$, the reverse statement in part 2 of Theorem \ref{thm:Mplus}
does not hold in general,
as the following example shows.
Let $M\in\RR^{m\times m}$, $m\ge2$, be a symmetric,
strictly positive definite matrix with identical diagonal elements,
$\tilde
\gamma\dvtx\RR^d\rightarrow\RR$ a variogram, and $C(x,y)=-M\tilde
\gamma
(x-y)$. Then
$C^{(0)}(x,y)$ given by (\ref{eq:x0}) is a cross covariance function, but
$\starexp(-M\tilde\gamma)$ is a positive definite function if and only
if
$\tilde\gamma\equiv0$. To see this, assume that $\starexp(-M\tilde
\gamma)$ is
a positive definite function and let $m=2$, $M=(M_{jk})_{j,k=1,2}$,
and $Z(x)=(Z_1(x), Z_2(x))$ be a corresponding
random field. Then with $a=(1, -1, 1, -1)^\top$ we have
\begin{eqnarray*}
\Var\bigl(Z_1(0) - Z_2(0) + Z_1(y) - Z_2(y)\bigr)
&=& a^\top\pmatrix{
\exp^*(-M\tilde\gamma(0)) & \exp^*(-M\tilde\gamma(y)) \cr
\exp^*(-M\tilde\gamma(y)) & \exp^*(-M\tilde\gamma(0))
} a
\\&=& 2(1, -1) \exp^*(-M\tilde\gamma(y)) (1, -1)^\top
\\&=& 4 \bigl(\mathrm{e}^{-M_{11} \tilde\gamma(y)} - \mathrm{e}^{-M_{12} \tilde\gamma(y)}\bigr)
.
\end{eqnarray*}
Since $M_{11} > M_{12}$,
the latter is non-negative if and only if $\tilde\gamma(y)=0$.

So, for an arbitrary cross variograms
$\gamma\dvtx\RR^{2d}\rightarrow\CC^{m\times m}$
the function $\starexp(-\gamma(x,y))$ is not a positive definite,
in general. However,
\[
C_1(x,y) = \starexp\bigl(\gamma(x, 0) + \gamma(y,0) - \gamma(x,y)\bigr)
\]
and
\begin{eqnarray}\label{eq:C2}
C_2(x,y) &=& \starexp\bigl(\gamma(x, 0) + \gamma(y,0) -D_{xy} - \gamma(x,y)\bigr),
\nonumber\\[-8pt]\\[-8pt]
(D_{xy})_{jk} &=& \gamma_{jj}(x,0) + \gamma_{kk}(y,0),\nonumber
\end{eqnarray}
are always positive definite functions in $\RR^d$,
cf. Theorem 2.2 in \citet{BCR} for the univariate case.
To see this, let $\gamma$ be an
$m$-variate cross variogram and $Z$ a corresponding
$m$-variate
random field. Let $Y(x) = Z(x) - Z(0)$ and $c(x,y) =\EE Y(x) Y^\top(y)$.
Then $c$ and $\overline{c^\top}$ are positive definite functions and
\begin{eqnarray*}
c_{jk}(x,y)+\overline{c_{kj}(x,y)}
&=& \EE\bigl(Y_j(x)\overline{Y_k(y)} +\overline{Y_k(x)}Y_j(y)\bigr)
\\&=&
\EE\bigl[ Y_j(x)\overline{Y_k(x)} + Y_j(y)\overline{Y_k(y)}
+ \bigl(Y_j(x) -Y_j(y)\bigr)\bigl(\overline{Y_k(y)}-\overline{Y_k(x)}\bigr)\bigr]
\\&=&\gamma_{jk}(x,0) + \gamma_{jk}(y,0) - \gamma_{jk}(x,y).
\end{eqnarray*}
Part 1 of Theorem 1 yields that $C_1$ is a positive definite function.
Let $Z$ be a corresponding random field. Then the random field
$(\mathrm{e}^{-\gamma_{11}(x,0)} Z_1(x),\ldots, \mathrm{e}^{-\gamma_{mm}(x,0)} Z_m(x))$,
$x\in\RR^d$, has cross covariance function $C_2$.
\end{remark}

\begin{remark}
Let $C(x_1,x_2)$ = $V D(x_1, x_2) \bar V^\top\in\CC^{m\times m}$,
$x,y\in\RR^d$,
for some unitary matrix $V\in\CC^{m\times m}$. The values of the mapping
$D\dvtx\RR^{2d}\rightarrow
\CC^{m\times m}$ are
diagonal matrices,
\[
D(x_1,x_2)=\diag(D_1(x_1, x_2),\ldots, D_m(x_1,x_2)),
\qquad x_1,x_2\in\RR^d,
\]
and the $D_j\dvtx\RR^{2d}\rightarrow\CC$, $j=1,\ldots, m$,
are arbitrary functions.
Then the $n$-fold matrix product
$C^n\dvtx\RR^{2d}\rightarrow\CC^{m\times m}$ is a
cross covariance function in $\RR^d$ for any $n\in\NN$
if and only if the $D_j$ are all covariance functions, and
Theorem~\ref{thm:Mplus} remains
true if $ \starexp(r C(x,y) )$ is replaced by
\[
\exp(r C(x,y) ) = \sum_{n=0}^\infty\frac{r^n C^{n}(x,y)}{n!} ,\qquad
x,y\in\RR^d.
\]
\end{remark}

The subsequent proposition generalizes
the results in \citet{cressiehuang99} and
Theorem~1 in \citet{gneitingnonseparable}.
Denote by $\B^d$ the ensemble of Borel sets of $\RR^d$.

\begin{proposition}\label{thm:bochner}
Let $d$ and $l$ be non-negative integers with $d+l>0$
and $C\dvtx\RR^{l + 2d}\rightarrow
\CC^{m\times m}$ a continuous function in the first argument.
Then the following two assertions are equivalent:
\begin{enumerate}
\item
$C$ is a cross covariance function that is translation invariant
in the first argument,
that is, $C(h, y_1, y_2) = \Cov(Z(x+h, y_1), Z(x, y_2))$ for some
second order random field $Z$ on $\RR^{l+d}$ and all
$x,h\in\RR^l$ and
$y_1,y_2\in\RR^d$.
\item
$C \dvtx \RR^l \times\RR^{2d}\rightarrow\CC^{m\times m}$
is the Fourier transform of some finite measures
$F_{y_1,y_2,j,k}$, $y_1,y_2 \in\RR^d$, $j,k=1,\ldots,m$,
that is,
\begin{equation}\label{eq:Cjk}
C_{jk}(h, y_1 , y_2) = \int \mathrm{e}^{-\mathrm{i} \langle h, \omega\rangle}
F_{y_1,y_2,j,k}(\DD\omega),\qquad h\in\RR^l, j,k=1,\ldots,m,
\end{equation}
and
\begin{equation}\label{eq:CA}
(C_{jk}^A(y_1,y_2))_{jk} = (F_{y_1,y_2,j,k}(A))_{jk},\qquad
y_1,y_2\in\RR^d,
\end{equation}
is an $m$-variate
cross covariance function in $\RR^d$ for any $A\in\B^l$.
\end{enumerate}
\end{proposition}

\begin{pf}
The proof follows the lines in \citet{gneitingnonseparable}.
Let us first assume that Equations~(\ref{eq:Cjk}) and (\ref{eq:CA}) hold.
Let $n\in\NN$, $x_1,\ldots, x_n\in\RR^l$, $y_1,\ldots, y_n\in\RR^d$
and $a_1,\ldots,
a_n\in\CC^m$ be fixed. Then a matrix-valued function $f\dvtx\RR^{l+
2d}\rightarrow
\CC^{m\times m}$
and a non-negative finite measure $F$ on $\RR^d$ exists, such
that
\begin{equation}\label{eq:fF}
\int_A f_{jk}(\omega, y_p, y_q)  F(\mathrm{d}\omega) = F_{y_p, y_q, j,
k}(A),\qquad
p,q=1,\ldots, n, j,k=1,\ldots,m,
\end{equation}
for any $A\in\B^l$. For instance, let
$F(A) = \sum_{p=1}^n \sum_{k=1}^m F_{y_p, y_p, k, k}(A)$.
Then, Equation~(\ref{eq:CA}) implies that the $mn\times mn$ matrix
$(f_{jk}(\omega, y_p, y_q))_{j,k; p,q}$
is hermitian for $F$-almost all $\omega$.
Now,
\begin{eqnarray*}
\sum_{p=1}^n \sum_{q=1}^n a_p^\top C(x_p - x_q, y_p, y_q) \overline{a_q}
&=&
\int\sum_{p=1}^n \sum_{q=1}^n \mathrm{e}^{ -\mathrm{i} \langle x_p,
\omega\rangle} a_p^\top
f(\omega, y_p, y_q)
\overline{ \mathrm{e}^{ -\mathrm{i} \langle x_q,\omega\rangle}} \overline{a_q}
F(\mathrm{d}\omega)
\ge0
.
\end{eqnarray*}

Conversely, let $C(h, y_1, y_2)\dvtx \RR^{l + 2d} \rightarrow\CC
^{m\times m}$
be a covariance function that is stationary in its first argument.
We have
\[
C_{jk}(h, y, y') = \int \mathrm{e}^{-\mathrm{i} \langle\omega, h\rangle}
F_{y, y', j, k}(\DD\omega),\qquad h \in\RR^l; y,y'\in\RR^d,
j,k=1,\ldots,m,
\]
for some finite, not necessarily positive measures $F_{y, y', j, k}$
(\citet{yaglomII}, page~115).
It now remains to demonstrate that equality (\ref{eq:CA}) holds.
Fix $n\in\NN$, $y_1,\ldots, y_n\in\RR^d$,
and $a_1,\ldots, a_n\in\CC^m$.
Then a non-negative finite measure $F$ and a function
$f\dvtx\RR^{l+ 2d}\rightarrow
\CC^{m\times m}$ exist, such that
Equation (\ref{eq:fF}) holds.
By assumption,
$
\sum_{p=1}^n \sum_{q=1}^n a_p^\top C(\cdot, y_p, y_q) \overline{a_q}
$
is a positive definite, continuous function and its
Fourier transform is non-negative.
Following directly from the linearity of the
Fourier transform, we have that for $F$-almost all $\omega\in\RR^l$
\[
\sum_{p=1}^n \sum_{q=1}^n a_p^\top f(\omega, y_p, y_q)
\overline{a_q}\ge0,
\]
which finally leads to Equation (\ref{eq:CA}).
\end{pf}

If a covariance function
is translation invariant, we
will write only one argument for ease of notation, for example,
$C(h)$, $h=x-y\in\RR^d$, instead of $C(x,y)$, $x,y\in\RR^d$.


\section{Generalized Gneiting's class}
\label{sec:gneiting}

A function
$C(x,y) = \varphi(\|h\|)$, $h=x-y\in\RR^d$,
is a motion invariant, real-valued covariance function in $\RR^d$ for
all $d\in\NN$ if and only if $\varphi$ is a normal scale mixture,
that is,
\[
\varphi(h) = \int_{[0,\infty)} \exp(- a h^2) \D F(a),\qquad
h\ge0,
\]
for some non-negative measure $F$ \citep{schoenberg38}.
Examples
are the stable model \citep{yaglomI},
the generalized Cauchy model \citep{cauchy},
\[
\varphi(h) = (1+h^\alpha)^{-\beta/\alpha},\qquad h\ge0
,
\]
$\alpha\in[0,2]$, $\beta>0$,
and the generalized hyperbolic model
(\citet{barndorffnielsen78}, \citet{gneitingdual}). The latter includes as
special case
the Whittle--Mat\'{e}rn model \citep{stein},
\[
\varphi(h) = W_\nu(h) = 2^{1-\nu} \Gamma(\nu)^{-1} h^\nu K_\nu
(h),\qquad h>0.
\]
Here, $\nu>0$ and $K_\nu$ is a modified Bessel function.

\begin{thm} \label{thm:main}
Assume that $m$ and $d$ are positive integers and
$H\dvtx\RR^d\rightarrow\RR^m$.
Suppose that $\varphi$ is a
normal scale mixture and $G\dvtx \RR^{2d} \rightarrow\RR^{m\times m}$ is
a cross variogram in $\RR^{d}$ or
$-G$ is a cross covariance function.
Let $M\in\RR^{m\times m}$ be positive definite,
such that $M + G(x, y)$
is strictly positive definite for all $x,y\in\RR^{d}$.
Then
\begin{equation}\label{eq:main}
C(x, y) = \frac{\varphi(
[(H(x)-H(y)) ^\top(M
+ G(x, y))^{-1}
(H(x)-H(y)) ]^{1/2}
)}
{\sqrt{|M + G(x, y)|}},\qquad
x,y\in\RR^d,
\end{equation}
is a covariance function in $\RR^d$.
\end{thm}

\begin{lemma}\label{lemma:sum}
Let $\gamma\dvtx \RR^{2d} \rightarrow\CC^{m\times m}$ be a
cross variogram (cross covariance function) in $\RR^d$
and $A\in\CC^{l \times m}$.
Then $\gamma_0 = A \gamma\overline{A^\top}$ is an
$l$-variate,
cross
variogram (cross covariance function)
in $\RR^{d}$.
\end{lemma}

\begin{pf*}{Proof of Theorem \ref{thm:main}}
We follow the proof in \citet{gneitingnonseparable} but assume
first that $\varphi(h) = \mathrm{e}^{-h^2}$.
If
$ G(x,y)$ is a cross variogram,
then, according to Lemma \ref{lemma:sum},
\[
g(x, y) = \omega^\top G( x, y) \omega
\]
is a (univariate) variogram for any $\omega\in\RR^m$.
Equation~(\ref{eq:C2}) or Theorem 2.2 in \citet{BCR} implies
\begin{equation}\label{eq:crucial.posdef}
C_\omega(x,y)= \exp(- \omega^\top G( x, y) \omega),\qquad
x,y\in\RR^d,
\end{equation}
and hence,
\begin{equation}\label{eq:crucial.posdef2}
\hat C(\omega, x,y)= \exp\bigl(- \omega^\top\bigl(M+G( x, y)\bigr)
\omega\bigr),\qquad x,y\in\RR^d,
\end{equation}
are both covariance functions for any fixed $\omega\in\RR^m$.
With $\D F_{x,y,1,1}(\omega) = \hat C(\omega, x,y) \D\omega$,
Proposition \ref{thm:bochner} yields that the univariate function
\begin{eqnarray*}
C(h, x, y)
& =&
c \frac{\exp(-
h^\top(M
+G(x,y))^{-1}
h
)}
{\sqrt{|M+ G( x, y)|}},\qquad
h \in\RR^m; x,y \in\RR^d
\end{eqnarray*}
is a covariance function in $\RR^{m+d}$ for all $c\ge0$, which
is translation invariant in the first argument.
Now, consider a random field $Z(\zeta, x)$
on $\RR^{m+d}$ corresponding to
$C(h, x, y)$ with $c=1$. Define the random field $Y$ on $\RR^d$ by
\[
Y(x) = Z(H(x), x)
.
\]
Then the covariance function of $Y$ is equal to the covariance function
given in the theorem.
For general $\varphi$, the assertion is obtained directly
from the definition of normal scale mixtures.
In case $-G$ is a cross covariance function, the proof runs exactly
the same way.
\end{pf*}

\begin{example}\label{ex:prod}
A well known construction of a cross covariance function in $\RR^d$
used in machine learning
is
\[
\tilde G(x,y) = f(x) f(y)^\top,\qquad x,y\in\RR^d,
\]
for some function $f\dvtx\RR^d \rightarrow\RR^{m \times l}$.
Assume that $M-f(x) f(y)^\top$ is strictly positive definite
for all $x$ and $y$ and some positive definite matrix $M$.
Then, $C$ in Equation~(\ref{eq:main}) is a covariance function
with $G = -\tilde G$.
\end{example}

We denote by $\1_{d\times d}\in\RR^{d\times d}$ the identity matrix.

\begin{example}
\citet{gneitingnonseparable}
delivers a rather general construction of non-separable models based
on completely monotone functions, containing as particular case
the models developed by
\citet{cressiehuang99}.
Let $\varphi$ be a completely monotone function, that is,
$ \varphi(t^2)$, $t\in\RR$,
is a normal scale mixture, and $\psi$ be
a positive function with a completely monotone derivative.
Then
\begin{equation}\label{eq:gneit}
C(h, u) = \frac1{\psi(|u|^2)^{d/2}} \varphi\bigl(\|h\|^2 /
\psi(|u|^2)\bigr),\qquad h \in\RR^d, u \in\RR,
\end{equation}
is a translation invariant covariance function in $\RR^{d+1}$
(\citet{gneitingnonseparable}, Theorem 2).
According to
Bernstein's theorem,
the function $\psi(\|\cdot\|^2)-c$
is a variogram for some positive constant $c$, see also \citet{BCR}.
The positive definite nature of $C$ in (\ref{eq:gneit}) is also
ensured by
Theorem~\ref{thm:main}
for $m=d$ and
$G((x_1, x_2), (y_1,y_2)) = \psi(\|x_2 -y_2\|^2) \1_{d\times d}$,
$x_1,y_1\in\RR^d$, $x_2,y_2\in\RR$.
\citet{gneitingnonseparable} provides examples for $\psi$ and,
along the way, introduces a new class of
variograms,
\[
\gamma(h) = (\|h\|^a + 1)^b -1,\qquad
a\in(0,2], b\in(0,1].
\]
This class generalizes the class of variograms of fractal Brownian motion
and that of multiquadric kernels \citep{wendland}.
\end{example}

\begin{example}\label{ex:cox}
In the context of modelling rainfall, \citet{coxisham88}
proposed in $\RR^{d+1}$ the translation invariant covariance function
\[
C(h, u) = \EE_V \varphi(\|h -V u\|),\qquad h \in\RR^d, u \in\RR
.
\]
Here, $\varphi(\|\cdot\|)$ is a motion invariant covariance function in
$\RR^d$
and $V$ is a $d$-dimensional random wind speed vector.
Unfortunately, this appealing model has lacked
explicit representations.
Now let us assume that $V$
\begin{figure}[b]

\includegraphics{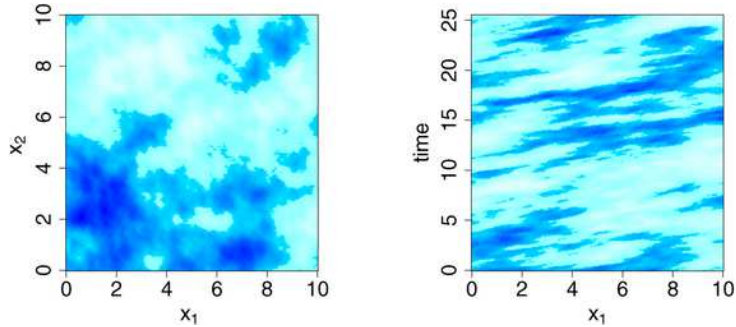}

\caption{Realizations of the Cox--Isham covariance model in $\RR
^2\times\RR$.
Left time $t=0$,
right $x_2=0$. See Example
\protect\ref{ex:cox} for details.}
\label{fig:rain}
\end{figure}%
follows a $d$-variate normal distribution $\N(\mu, D/2)$
and $\varphi(x)=\exp(-x^2)$.
Then,
\[
C(h, u) = \frac1{\sqrt{|\1_{d\times d}+ u^2 D|}}
\varphi\bigl(
[ (h - u \mu)^\top(\1_{d\times d} + u^2 D)^{-1} (h-u\mu) ]^{1/2}
\bigr),\qquad
h\in\RR^d, u\in\RR
,
\]
please refer to the appendix for a proof.
Hence, $C(h,u)$ above is a covariance function
for any normal mixture $\varphi$.
Figure \ref{fig:rain} provides realizations of a random field
with the above covariance function
where $\varphi=W_1$ is the Whittle--Mat\'{e}rn model,
$\mu=(1,1)$ and
\[
D= \pmatrix{
1 & 0.5 \cr
0.5 & 1
}.
\]
\end{example}

\begin{remark}\label{ex:stein1}
\citet{stein05b} proposes models in $\RR^d$ of the form
\[
C(x, y)
= \frac{\varphi(
[(x-y)^\top(f(x) + f(y))^{-1}
(x-y) ]^{1/2}
)}
{\sqrt{|f(x) + f(y)|}},\qquad
x,y\in\RR^d,
\]
in which the values of $f\dvtx \RR^{2d} \rightarrow\RR^{m\times m}$ are
strictly positive definite matrices, see also \citet{paciorek03}
and \citet{PMC09}.
Here, $f(x) + f(y)$ is not
a variogram in general, but the proof of Theorem \ref{thm:main}
is still applicable
if $\hat C$ in Equation~(\ref{eq:crucial.posdef2}) is replaced by
\[
\hat C(\omega, x,y) = \exp\bigl(- \omega^\top\bigl(f(x) + f(y)\bigr) \omega\bigr)
,
\]
which is a positive definite function for all $\omega\in\RR^m$.
\end{remark}

\begin{remark}
The covariance model (\ref{eq:main}), which is valid in $\RR^d$,
does not allow for negative values, hence
its value is limited in some applications \citep{GPMS07}.
To overcome this limitation, \citet{ma05}
considers differences of positive definite
functions.
Let $B_1,B_2,M_1,M_2\in\RR^{d\times d} $
be strictly positive definite matrices.
Proposition \ref{thm:bochner} yields that
\begin{eqnarray*}
C(h, x, y) &=& \frac{\exp(- [h^\top(M_1 + (x-y)^\top B_1
(x-y)\1_{d\times d})^{-1}
h ])}{\sqrt{|M_1 + (x-y)^\top B_1 (x-y)\1_{d\times d}|}}
\\&&{} + b\frac{\exp(- [h^\top
(M_2 + (x-y)^\top B_2 (x-y)\1_{d\times d})^{-1}
h ])}{\sqrt{|M_2 + (x-y)^\top B_2 (x-y)\1_{d\times d}|}} ,\qquad
h,x,y\in\RR^d,
\end{eqnarray*}
is a positive definite function in $\RR^{2d}$ that is translation
invariant in
its first argument
if and only if for all $\omega\in\RR^d$,
\begin{eqnarray*}
\hat C_\omega(x, y) &=& \exp\bigl(- \omega^\top M_1\omega
- \|\omega\|^2 (x-y)^\top B_1
(x-y)
\bigr)
\\&&{} + b\exp\bigl(- \omega^\top M_2\omega
- \|\omega\|^2 (x-y)^\top B_2
(x-y)
\bigr),\qquad
x,y\in\RR^d,
\end{eqnarray*}
is a positive definite function, that is, if and only if for all
$\omega
,\xi\in\RR^d$,
\begin{eqnarray*}
|B_1|^{-1/2} \exp(- \omega^\top M_1\omega
- \|\omega\|^2 \xi^\top B_1^{-1}
\xi
)
+ b |B_2|^{-1/2}\exp(- \omega^\top M_2\omega
- \|\omega\|^2 \xi^\top B_2^{-1}
\xi
) \ge0.
\end{eqnarray*}
This is true for some negative value of $b$ if and only if
both $M_2 -M_1$ and $B_2^{-1}-B_1^{-1}$ are positive definite matrices.
In this case, $C(h,x,y)\dvtx\RR^{3d}\rightarrow\RR$ is
a positive definite function in $\RR^{2d}$
if and only if
\[
b \ge- \sqrt{ |B_2| / |B_1|}
.
\]
Then, $C_0$ given by
$ C_0(x,y) = C(x-y,x, y) $
is a stationary
covariance function in $\RR^d$ that may take negative values.
\end{remark}

\begin{remark}
The condition that $M+G(x,y)$ is strictly positive definite
for all $x,y\in\RR^d$ can be relaxed.
For example, let $d=2$ and $(h,u) =x-y\in\RR^2$. Then, the function
$C(h,u) = |u|^{-1/2} \exp(-h^2 / |u|)$
is of the form (\ref{eq:main}) and
defines a covariance function of
a stationary, generalized random field on $\RR^2$, see Chapter 3 in
\citet{gelvandvilenkin4} and Chapter 17 in \citet{koralovsinai}.
Note that, here, $\lim_{u\rightarrow0} C(0,u) = \infty$. Hence, $C$ cannot
be a translation invariant
covariance function in the usual sense.
\end{remark}


\section{Model constructions based on dependent processes}
\label{sec:random}

The idea of the subsequent two constructions is based on the following
observation. Let $C(h,u)=C_0(h) C_1(u)$, $h\in\RR^d$, $u\in\RR$,
be a translation invariant, real-valued covariance model
in $\RR^{d+1}$
and assume we are interested in the corresponding random field
at some fixed locations $x_1,\ldots, x_n\in\RR^d$
and for all $t\in\RR$.
Let $Y_x$, $x\in\RR^d$, be i.i.d.\ temporal processes with
covariance function $C_1$. Then
\[
Z(t) =(Z_{x_1}(t),\ldots,Z_{x_n}(t))
= \bigl(C_0(x_p-x_q)\bigr)_{p,q=1,\ldots,n}^{1/2} (Y_{x_1}(t), \ldots,
Y_{x_n}(t))^\top,\qquad
t\in\RR,
\]
has the required covariance structure.
Now, $Z$ can be interpreted as a finite, weighted sum
over $Y_x$, $x\in\RR^d$.
The separability is caused by the fact that $Y$ enters into
the sum only through the fixed instance $t$.
Non-separable models can be obtained if the argument of $Y$
also depends on the location.

\subsection{Moving averages based on fields of temporal processes}
Assume that
$Y(A, t)$, $A\in\B^d$ and $t\in\RR^l$,
is a stationary process such that $Y(A_1, \cdot),\ldots,$ $Y(A_n,
\cdot)$
are independent for any disjoint sets $A_1,\ldots,A_n\in\B^d$, $n\in
\NN
$. In the second argument, $Y$ is
a stationary, zero mean Gaussian random field on $\RR^{l}$
with covariance function $|A|C_1$, $C_1\dvtx\RR^l \rightarrow\RR$.
Then,
\[
\Cov(Y(A,t), Y(B,s)) = |A \cap B| C_1(t -s)
\]
for any $s,t \in\RR^l$ and $ A, B \in\B^d$.
Let $f\dvtx\RR^d\rightarrow\RR^l$ be continuous,
$g\dvtx\RR^d\rightarrow\RR$ be continuous and
square-integrable, and
\[
Z(x,t) = \int_{\RR^d} g(v-x) Y\bigl(\D v, f(v-x) - t\bigr),\qquad
x\in\RR^d, t \in\RR^l.
\]
Then $Z$ is weakly stationary on $\RR^{d+l}$ with translation invariant
covariance function
\begin{eqnarray*}
C(h, u) &=&
\int_{\RR^d} g(v) g(v+h) C_1\bigl(f(v) - f(v+h) -u\bigr) \D v,\qquad h\in\RR^d, u \in\RR^l.
\end{eqnarray*}

\begin{example}\label{ex:fields}
Let $g(v) = (2\uppi^{-1})^{d/4} \exp(- \|v\|^2)$, $v\in\RR^d$,
$l=1$, $C_1(u) = \exp(-u^2)$, $u\in\RR$, and
$f(v)= v^\top A v +z^\top v$, $v\in\RR^d$,
for a symmetric, not necessarily positive
definite matrix $A \in\RR^{d\times d}$ and $z\in\RR^d$.
Let us further introduce a non-negative random scale $V$, that is,
\[
Z(x,t) =
V^{d/2}
\int_{\RR^d} g\bigl(\sqrt{V}(v-x)\bigr) Y\bigl( \DD v, \sqrt{V}\bigl(f(v-x) - t\bigr)\bigr),\qquad
x\in\RR^d, t \in\RR.
\]
Let $B=A h h^\top A$.
Then the covariance function of $Z$ equals
\begin{eqnarray}\label{eq:fields}
C(h, u) =
|\1_{d\times d} + 2 B|^{-1/2}
\EE_V \mathrm{e}^{-V[
\|h\|^2/2 + (z^\top h+ u)^2 (1 - 2 h^\top A(\1_{d\times d} + 2
B)^{-1}A h)
]}
,
\end{eqnarray}
please refer to the appendix for a proof.
Equation~(\ref{eq:fields}) reveals that $C$ is a potential covariance model
for rainfall with frozen wind direction.
Figure \ref{fig:moving} depicts realizations
\begin{figure}

\includegraphics{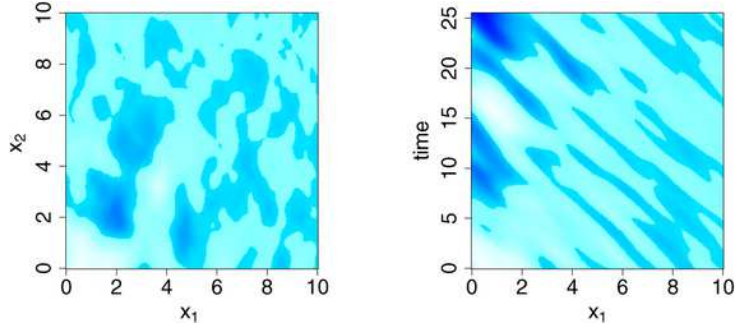}

\caption{Realizations of a moving average random field in $\RR
^2\times
\RR$.
Left time $t=0$, right $x_2=0$.
See Example \protect\ref{ex:fields} for the definition of the covariance structure.}
\label{fig:moving}
\end{figure}%
of a random field with the above covariance function
where
$\EE_V \exp(-V Q)$ is the Whittle--Mat\'{e}rn model $W_1(\sqrt Q)$,
$Q\ge0$,
$z=(2,0)$ and
$
A = \left({0.5 \atop 0}\enskip {0 \atop 1}\right)$.
\end{example}

\subsection{Models based on a single temporal process}
Another class of models may be
obtained by considering only a single process $Y$.
Although the subsequent approach might be generalized,
an explicit model has currently only been found
within the framework of normal scale mixtures.
For $x \in\RR^d$ let
\begin{equation}\label{eq:single}
Z(x) = (2V/\uppi)^{d/4}|S_x|^{1/4}
\mathrm{e}^{- V (U-x)^\top S_x (U-x)}
Y \bigl(\sqrt V \bigl(\xi_1(U-x) + \xi_2(x) \bigr) \bigr)
\frac{g(V, x) }{ \sqrt{f(U)}}.
\end{equation}
Here, $V$ is a positive random variable and $U$ is
a $d$-dimensional random variable with strictly positive density $f$.
The one-dimensional random process $Y$ is assumed to be stationary
with Gaussian covariance function $C(t) = \mathrm{e}^{-t^2}$. The matrix $S_x$ is
strictly positive definite
for all $x\in\RR^d$,
$ \xi_2\dvtx\RR^d \rightarrow\RR$
is arbitrary, and
$g$ is a positive function such that $\EE_V g(V,x)^2$ is finite for all
$x\in\RR^d$.
The function $\xi_1$ is quadratic, that is,
\[
\xi_1(x)= x^\top M x+ z^\top x
\]
for a symmetric $d\times d$ matrix $M$ and an arbitrary vector $z\in
\RR^d$.
Let
\begin{eqnarray*}
c &=& -z^\top(x-y) +\xi_2(x) - \xi_2(y),
\\
A&=&S_x + S_y + 4M (x-y) (x-y) ^\top M,
\\
m &=&(x-y)^\top M (x-y),
\end{eqnarray*}
and
\[
Q(x,y) =  c^2- m^2
+ (x-y)^\top\bigl(S_x + 2 (m +c) M\bigr) A^{-1}
\bigl(S_y + 2 (m - c) M\bigr) (x-y).
\]
Then the covariance function of $Z$ equals
\begin{equation}\label{eq:single2}
C(x,y) = \frac{2^{d/2}
|S_x|^{1/4}|S_y|^{1/4}}{\sqrt{|A|}}\cdot
\EE_V g(V, x)g(V, y)\exp(-V Q(x,y) )
,\qquad x,y\in\RR^d.
\end{equation}
The proof is given in the \hyperref[app]{Appendix}.

\begin{example}\label{ex:single}
Translation-invariant models in $\RR^d$ are obtained
if both $S_x$ and $g$ do not depend on~$x$.
Assume $S_x$ is twice the identity matrix,
$g(v) = (2\sqrt{v})^{1-\nu} / \sqrt{\Gamma(\nu)}$, $v,\nu>0$,
and $V$ follows the Fr\'{e}chet distribution
$F(v) = \mathrm{e}^{-1/(4v)}$, $v>0$.
Two particular models might be of special interest, either because of
their
simplicity or their explicit spatio-temporal modelling.
First, if $c\equiv0$ then
\[
C(h) = \frac{W_\nu(\|h\|)}{|\1_{d\times d} + M h h^\top M|^{1/2}},
\qquad
h\in\RR^d,
\]
according to formula \textup{3.471.9} in Gradshteyn and Ryzhik (\citeyear{GRengl}).
Second, an explicit spatio-temporal model in $\RR^{d+1}$ is obtained
for
\[
\xi_2(x,t) = t,\qquad x\in\RR^d, t\in\RR,
\quad\mbox{and}\quad M = \pmatrix{
L &  0\cr
0 & 0
}.
\]
Then, with $D=\1_{d\times d} + L h h^\top L$, we get
\[
C(h, u)
= |D|^{-1/2} W_\nu\bigl(\sqrt{Q(h,u)}\bigr),\qquad
h\in\RR^d, u\in\RR,
\]
where
\[
Q(h,u) = (u -z^\top h)^2 - (h^\top L h)^2
+ h^\top\bigl(D + (u -z^\top h)L\bigr) D^{-1}
\bigl(D + (u -z^\top h)L\bigr)h.
\]
\end{example}

\begin{example}
Let $\xi_1\equiv\xi_2\equiv0$.
Then the random process $Y(t)$ is considered only at instance $t=0$ and
the exponent $Q(x,y)$ simplifies to
\[
Q(x, y) = (x-y)^\top S_x (S_x + S_y)^{-1} S_y (x-y)
= (x-y)^\top(S_x^{-1} + S_y^{-1} )^{-1}(x-y)
.
\]
Let $g(v,x) = (2\sqrt v)^{1-\nu(x)} / \Gamma(\nu(x))^{1/2}$,
$\nu$ a positive function on $\RR^d$, and $V$
a Fr\'{e}chet variable with distribution function $F(v) = \mathrm{e}^{-1/(4 v)}$,
$v>0$.
Then, the first model given in \citet{stein05b} is obtained,
\[
C(x, y) = \frac{2^{d/2}
|S_x|^{1/4}|S_y|^{1/4} \Gamma((\nu(x) + \nu(y)) /2)}
{[|S_x + S_y|\Gamma(\nu(x)) \Gamma(\nu(y))]^{1/2}}
W_{(\nu(x) + \nu(y)) /2}(Q(x,y)^{1/2}),\qquad
x,y\in\RR^d
.
\]
The second model given in \citet{stein05b},
a generalization of the Cauchy model, is obtained by
$g(v, x) = v^{(\delta(x) -1) / 2}$ and a standard exponential
random variable $V$, that is,
\[
C(x, y) = \frac{2^{d/2}
|S_x|^{1/4}|S_y|^{1/4}}
{|S_x + S_y|^{1/2}(1 + Q(x,y))^{(\delta(x) + \delta(y))/2}},\qquad
x,y\in\RR^d.
\]

If $\nu$ and $\delta$ are constant, then the above models are special cases
of Theorem
\ref{thm:main}.

See Theorem 1 in \citet{PMC09} for a class of models that generalizes
Stein's examples.
\end{example}

\begin{example}\label{ex:cyclone}
A cyclone can be mimicked if rotation matrices are included in
the model,
\[
C(x,y) = \frac{2^{d/2}|S_x|^{1/4}|S_y|^{1/4}}{\sqrt{|S_x + S_y|}}
W_\nu\bigl(
\bigl( h^\top S_x (S_x + S_y)^{-1} S_y h \bigr)^{1/2}
\bigr),\qquad
x,y,\in\RR^3,
\]
where
\begin{eqnarray*}
S_x&= &\diag( 1, 1, 1) + R(x)^\top A^\top x x^\top A R(x),\qquad A \in
\RR^{3 \times3},
\\
R(x)&=& \pmatrix{
\cos(\alpha x_3) & -\sin(\alpha x_3) & 0
\cr
\sin(\alpha x_3) & \cos(\alpha x_3)& 0
\cr
0& 0& 1
},\qquad x = (x_1, x_2, x_3)\in\RR^3, \alpha\in\RR,
\end{eqnarray*}
and
\[
h  =  x^\top R(x) - y^\top R(y).
\]
The positive definiteness
of the model is now ensured by both Theorem
\ref{thm:main} and a generalized version of $Z$ in Equation~(\ref{eq:single}),
replacing $x$ by $x^\top R(x)$ there.
Note that $x\mapsto x^\top R(x)$ is a bijection.
\begin{figure}

\includegraphics{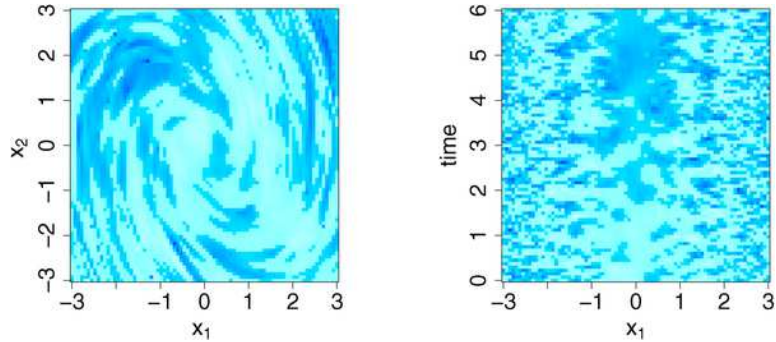}

\caption{Realizations of a random field in $\RR^3$
that mimics a cyclone. Left time $x_3=0$, right $x_2=0$.
See Example \protect\ref{ex:cyclone}
for the definition of the covariance structure.}
\label{fig:hurricane}
\end{figure}%
Figure \ref{fig:hurricane} depicts realizations of a random field
with the above covariance function where
$\alpha= -2 \uppi$, $ \nu=1$, and
\[
A= \pmatrix{
2 & 1 & 0\cr
0 & 1 & 0\cr
0 & 0 & 0
}.
\]
\end{example}


\section{Multivariate spatio-temporal models}
\label{sec:multivariate}
Here, we generalize Theorem \ref{thm:main} to construct multivariate cross
covariance functions.
Let $\underline M = (M + M^\top) /2$ for any real-valued square matrix~$M$.

\begin{thm} \label{thm:multi}
Assume that $l$, $m$ and $d$ are positive integers,
$A_{j}\in\RR^{l\times d}$ for $j=1,\ldots,m$.
Suppose that $\varphi$ is a
normal scale mixture and $G\dvtx \RR^{2 d} \rightarrow\RR^{l\times l}$
is a cross covariance function.
Let $M\in\RR^{d\times d}$ be a positive definite matrix
such that $M - \underline{A_j^\top G(x, y)A_k}$
is strictly positive definite for all $x,y\in\RR^{d}$ and
$j,k=1,\ldots,d$.
Then $C=(C_{jk})_{j,k=1,\ldots,m}$ is a
cross covariance function in $\RR^d$ for
\begin{eqnarray}\label{eq:extension1}
&&C_{jk}(x, y)=
\frac{\varphi(
[(x-y) ^\top
(M
- \underline{A^\top_{j}G(x, y)
A_{k}})^{-1} (x-y) ]^{1/2}
)}
{\sqrt{|M - \underline{A^\top_{j}G(x, y)
A_{k}}|}},\nonumber\\[-8pt]\\[-8pt]
&&\quad
x,y\in\RR^d, j,k=1,\ldots,m.\nonumber
\end{eqnarray}
\end{thm}

\begin{pf}
Lemma \ref{lemma:sum} yields that
\begin{eqnarray*}
(\omega^\top\underline{A_j^\top G(x,y) A_k} \omega)_{j,k=1,\ldots,m}
&=&
(\omega^\top A_j^\top G(x,y) A_k \omega)_{j,k=1,\ldots,m}
\\&= &
(A_1\omega, \ldots, A_m\omega)^\top G(x,y) (A_1\omega, \ldots,
A_m\omega)
\end{eqnarray*}
is a cross covariance function for all $\omega\in\RR^d$.
Part 1 of Theorem \ref{thm:Mplus} yields that
$C_\omega(x,y)
\kern-0.5pt = \kern-0.5pt
(\exp(\omega^\top\kern-2pt\underline{A_j^\top G(x,y) A_k}
\omega))_{j,k=1,\ldots,m}$
is also a cross covariance function.
By assumption, $M- \underline{A_j^\top G(x,y) A_k} $ is strictly
positive definite.
Hence, as a result of Proposition \ref{thm:bochner},
the Fourier transform of the function $\omega\mapsto\exp(-\omega
^\top M
\omega)C_\omega(x,y)$ is a cross covariance
function, which is of the form (\ref{eq:extension1}).
\end{pf}


\begin{appendix}
\section*{Appendix}\label{app}
\renewcommand{\theequation}{\arabic{equation}}

\subsection{\texorpdfstring{Proof for the covariance function in Example
\protect\ref{ex:cox}}{Proof for the covariance function in Example 9}}
Let $f_{\mu, D/2} (x)$ be the multivariate normal density
with expectation $\mu$ and covariance matrix $D/2$.
Then we get
\begin{eqnarray*}
&&- \log\bigl( \varphi(h-uv)f_{\mu, D/2} (v)\bigr) + \tfrac12 \log
((2\uppi)^d |D|)
\\
&&\quad=
h^\top h - 2u h^\top v + u^2 v^\top v + v^\top D^{-1} v - 2\mu^\top
D^{-1} v +
\mu^\top D^{-1} \mu
\\
&&\quad=
h^\top h + \mu^\top D^{-1} \mu+ (v-\xi)^\top
(u^2 \1_{d\times d} + D^{-1})(v-\xi) -\xi^\top(u^2 \1_{d\times d} +
D^{-1})\xi
\end{eqnarray*}
with $\xi= (u^2 \1_{d\times d} + D^{-1})^{-1}(u h + D^{-1} \mu)$.
Hence,
\begin{eqnarray*}
&&
-\log C(h,u) + \tfrac12 \log(|D|) +\tfrac12\log(|u^2 \1_{d\times d} + D^{-1}|)
\\
&&\quad=
h^\top h + \mu^\top D^{-1} \mu-\xi^\top(u^2 \1_{d\times d} +
D^{-1})\xi
\\
&&\quad=
(h- u \mu)^\top(\1_{d\times d} + u^2 D)^{-1} (h- u \mu)
\end{eqnarray*}
which yields the assertion.

\subsection{\texorpdfstring{Proof for the covariance function in Example
\protect\ref{ex:fields}}{Proof for the covariance function in Example 13}}
We proof the formula for the covariance function in Example \ref{ex:fields},
but also demonstrate that a slightly more general function $g$ does
not give a more general
model.
To this end, let $g(v) = (|2\uppi^{-1}M|)^{1/4} \exp(- v^\top M v)$,
$v\in\RR^d$,
for a strictly
positive definite matrix $M \in\RR^{d\times d}$.
For ease of notation we assume that $V\equiv1$.
Then
\begin{eqnarray*}
&&-\log\bigl(g(v) g(v+h) C_1\bigl(f(v) - f(v+h) -u\bigr)\bigr)
- \tfrac12 \log(|2\uppi^{-1}M|)
\\
&&\quad=
v^\top M v + (v+h)^\top M (v + h) +
(2 v^\top A h + h^\top A h +z^\top h +u)^2
\\
&&\quad=
2 v^\top M v + 4 v^\top B v + 2 v^\top(2B + M +2uA + 2Ahz^\top)h
+ c
\end{eqnarray*}
where $B = Ahh^\top A$ and
$c = [h^\top Ah +z^\top h+ u ]^2 + h^\top M h$.
Hence, with $D = 2B + M +2[u + z^\top h]A$,
\begin{eqnarray*}
&&-\log\bigl(g(v) g(v+h) C_1\bigl(f(v) - f(v+h) + u\bigr)\bigr)
- \tfrac12 \log(|2\uppi^{-1}M|)
\\
&&\quad=
\bigl(v - (2M+4B)^{-1} Dh\bigr)^\top(2M + 4B) \bigl(v - (2M+4B)^{-1}
Dh\bigr)\\
&&\qquad
{}- h^\top D (2M + 4B)^{-1} D h + c.
\end{eqnarray*}
Thus,
\begin{eqnarray*}
C(h, u) &=& \frac{|M|^{1/2}}{|M + 2 B|^{1/2}}
\exp\bigl(-c + h^\top D (2M + 4B)^{-1} D h\bigr),\qquad
h\in\RR^d, u\in\RR
.
\end{eqnarray*}

Let $M^{-1/2}$ be a symmetric matrix with $M^{-1/2} M M^{-1/2}=\1
_{d\times d}$.
Replacing on the right hand side
$M^{-1/2} A M^{-1/2}$ by $\tilde A$, $M^{-1/2}z$ by $\tilde z$ and $M^{1/2}
h$ by $\tilde h$
shows that $M$ causes nothing but a geometrical anisotropy effect.
Hence, we may assume that
$M$ is the identity matrix. Then
\[
C(h, u) = |\1_{d\times d} + 2 B|^{-1/2}
\exp\bigl(-\bigl[c - \tfrac12h^\top D (\1_{d\times d} + 2B)^{-1} D h\bigr]\bigr)
\]
which yields Equation~(\ref{eq:fields}).

\subsection{\texorpdfstring{Proof of
Equation~(\protect\ref{eq:single2})}{Proof of Equation~(13)}}

Let $h=x-y$ and $w = U-x$. Then we have
\begin{eqnarray*}
\Cov(Z(x), Z(y)) &=&
\uppi^{-d/2}|S_x|^{1/4}|S_y|^{1/4}\EE_V V^{d/2}g(V, x)g(V, y)
\\
&&
{}\times
\int
\exp\bigl(-V w^\top S_x w - V (w+h)^\top S_y(w+h)\\
&&
\hphantom{{}\times
\int
\exp\bigl(}
{}- V \bigl( w^\top M w - (w+h)^\top M (w+h) + c \bigr)^2
\bigr)\D w.
\end{eqnarray*}
The value of the integral is at most
$\int\exp(-V w^\top S_x w) \D w$. Hence $\Cov(Z(x), Z(y))<\infty$
if $\EE_V g(V, x)g(V, y) < \infty$.
Now,
\begin{eqnarray*}
&&w^\top S_x w + (w+h)^\top S_y(w+h) + \bigl(w^\top M w -
(w+h)^\top
M (w+h)
+c\bigr)^2
\\
&&\quad=
w^\top(S_x + S_y + 4 M h h^\top M)w +
2w^\top\bigl(S_y + 2 (h^\top M h -c) M\bigr) h + h^\top S_y h +(h^\top M h -c)^2
\\
&&\quad=
(w - \mu)^\top A (w- \mu) -\mu^\top A \mu+ h^\top S_y h +(h^\top M
h -c)^2
\end{eqnarray*}
with $\mu= -A^{-1}(S_y + 2 (h^\top M h -c) M)h$.
That is,
\begin{eqnarray}
\label{eq:13A}
\Cov(Z(x) ,Z(y)) &=& |A|^{-1/2}{ |S_x|^{1/4}|S_y|^{1/4}}{}\EE_V
g(V, x)g(V, y)\nonumber\\[-8pt]\\[-8pt]
&&{}\times
\mathrm{e}^{-V [ hS_y h + (h^\top M h -c)^2 - \mu^\top A \mu]}.\nonumber
\end{eqnarray}
On the other hand, using the transform $w=U-y$, we get
\begin{eqnarray}\label{eq:13B}
&&\Cov(Z(x), Z(y))\nonumber
\\
&&\quad=
\uppi^{-d/2}|S_x|^{1/4}|S_y|^{1/4}\EE_V V^{d/2} g(V, x)g(V,
y)\nonumber\\
&&\qquad
{}\times\int
\exp\bigl(-V(w-h)^\top S_x (w-h) +
- V hS_y h\\
&&
\qquad\hphantom{{}\times\int
\exp\bigl(}
{}- V \bigl( (w-h)^\top M (w-h) - w^\top M w + c \bigr)^2
\bigr)
\D w\nonumber
\\&&\quad=
|A|^{-1/2} |S_x|^{1/4}|S_y|^{1/4}
\EE_V g(V, x)g(V, y)
\mathrm{e}^{-V [ hS_x h + (h^\top M h +c)^2 - \nu^\top A \nu]}\nonumber
\end{eqnarray}
with $\nu=A^{-1}(S_x + 2 (h^\top M h +c) M)h$.

Choosing $V\equiv1$ and $g$ a constant function we obtain that
the exponents in (\ref{eq:13A}) and (\ref{eq:13B}) must be equal, that
is,
\begin{eqnarray*}
&&hS_y h + (h^\top M h -c)^2 - \mu^\top A \mu
\\
&&\quad=
\tfrac12 [ hS_y h + (h^\top M h -c)^2 - \mu^\top A \mu
+ hS_x h + (h^\top M h +c)^2 - \nu^\top A \nu]
\\
&&\quad=
\tfrac12 [h(S_y + S_x + 4Mhh^\top M) h - 2 (h^\top M h)^2 + 2 c^2
- (\mu- \nu) A (\mu- \nu) - 2 \nu^\top A^{-1} \mu]
\\
&&\quad=
c^2- (h^\top M h)^2 - \nu^\top A^{-1} \mu.
\end{eqnarray*}
\end{appendix}

\section*{Acknowledgements}
The author is grateful to Zakhar Kabluchko, Emilio Porcu
and the referees for
valuable suggestions and comments.

\printhistory

\end{document}